\documentclass[letterpaper, 10 pt, conference]{ieeeconf}  

\IEEEoverridecommandlockouts                              

\usepackage{amsfonts}
\usepackage{graphics} 
\usepackage{epsfig} 
\usepackage{amsmath} 
\usepackage{amssymb}  
\usepackage{tikz}
\usepackage{tikz-cd}

\newtheorem{definition}{Definition}
\newtheorem{prob}{Problem}

\title{\LARGE \bf
Optimality of Zeno Executions in Hybrid Systems*
}

\author{William Clark$^{1}$ and Maria Oprea$^{2}$
\thanks{*W. Clark was funded by NSF grant DMS-1645643. M. Oprea was supported by the Army Research Office Biomathematics Program Grant W911NF-18-1-0351.}
\thanks{$^{1}$William Clark is with the Department of Mathematics,
        Cornell University, Ithaca, NY, 14580, USA
        {\tt\small wac76@cornell.edu}}%
\thanks{$^{2}$Maria Oprea is with the Center of Applied Mathematics, Cornell University, Ithaca, NY, 14580, USA
        {\tt\small mao237@cornell.edu}}
}

\begin{document}

\maketitle
\thispagestyle{empty}
\pagestyle{empty}

\begin{abstract}

A unique feature of hybrid dynamical systems (systems whose evolution is subject to both continuous- and discrete-time laws) is Zeno trajectories. Usually these trajectories are avoided as they can cause incorrect numerical results as the problem becomes ill-conditioned. However, these are difficult to justifiably avoid as determining when and where they occur is a non-trivial task. It turns out that in optimal control problems, not only can they not be avoided, but are sometimes required in synthesizing the solutions. This work explores the pedagogical example of the bouncing ball to demonstrate the importance of ``Zeno control executions.''

\end{abstract}

\section{INTRODUCTION}
A hybrid system is a system whose dynamics are controlled by a mixture of both continuous and discrete transitions. For the purposes of this note, such a dynamical system will be described via
\begin{equation}\label{eq:HDS}
	\begin{cases}
		\dot{x} = X(x), & x\in M\setminus \mathcal{S}, \\
		x^+ = \Delta(x^-), & x\in\mathcal{S}.
	\end{cases}
\end{equation}
The set $\mathcal{S}$ will be referred to as the guard and $\Delta$ as the reset. All of the data will be tacitly assumed to be sufficiently smooth. A unique phenomenon to hybrid systems is that of \textit{Zeno}. A trajectory, $\gamma(t)$, is Zeno if it intersects the guard infinitely many times in a finite amount of time. Determining when/where this occurs is a difficult problem \cite{1582237,4434891,10.1007} and is commonly excluded. One way of ensuring that Zeno does not occur is by requiring that $\overline{\mathcal{S}}$ and $\overline{\Delta(\mathcal{S})}$ are disjoint and the set of impact times is closed and discrete \cite{9838077,8836630,898695}.

In the case of mechanical systems with impacts, the state-space is the tangent bundle of the configuration space, $M=TQ$, and the guard consists of all outward pointing vectors at the location of impact, i.e.
\begin{equation*}
	\mathcal{S} = \left\{ (x,v)\in TQ : h(x)=0, \; dh_x(v) < 0 \right\},
\end{equation*}
where the condition $h(x)=0$ describes the location of impacts. As mechanical impacts reverse the normal component of velocity, 
\begin{equation*}
	\Delta(\mathcal{S}) = \left\{ (x,v)\in TQ : h(x)=0,\; dh_x(v)>0 \right\}.
\end{equation*}
In this setting, $\overline{\mathcal{S}}$ and $\overline{\Delta(\mathcal{S})}$ are \textit{never} disjoint and Zeno remains a possibility.  Although Zeno behaviour  will almost never occur for elastic impacts \cite{2101.11128},  for inelastic mechanical impacts, Zeno behavior is to be expected.

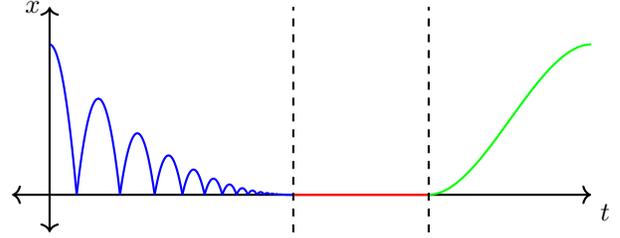
\begin{figure}
	\centering
	\begin{tikzpicture}[xscale=0.36]
			\draw[<->,thick] (-1.39,0) -- (20,0); 
			\draw[<->,thick] (0,-0.5) -- (0,2.5);
			\node [below right] at (20,0) {$t$};
			\node [left] at (0,2.5) {$x$};
			\draw[thick, blue, domain=0:1] plot (\x,{-2*\x*\x+2});
			\draw[thick, blue, domain=0:1.6] plot (\x+1,{-2*\x*\x+3.2*\x});
			\draw[thick, blue, domain=0:1.28] plot (\x+2.6,{-2*\x*\x+2.56*\x});
			\draw[thick, blue, domain=0:1.024] plot (\x+3.88,{-2*\x*\x+2.048*\x});
			\draw[thick, blue, domain=0:0.8192] plot (\x+4.904,{-2*\x*\x+1.6384*\x});
			\draw[thick, blue, domain=0:0.6554] plot (\x+5.7232,{-2*\x*\x+1.3107*\x});
			\draw[thick, blue, domain=0:0.5243] plot (\x+6.3786,{-2*\x*\x+1.0486*\x});
			\draw[thick, blue, domain=0:0.4194] plot (\x+6.9029,{-2*\x*\x+0.8389*\x});
			\draw[thick, blue, domain=0:0.3356] plot (\x+7.3223,{-2*\x*\x+0.6711*\x});
			\draw[thick, blue, domain=0:0.2684] plot (\x+7.6579,{-2*\x*\x+0.5369*\x});
			\draw[thick, blue, domain=0:0.2148] plot (\x+7.9263,{-2*\x*\x+0.4295*\x});
			\draw[thick, blue, domain=0:0.1718] plot (\x+8.1411,{-2*\x*\x+0.3436*\x});
			\draw[thick, blue, domain=0:0.1374] plot (\x+8.3129,{-2*\x*\x+0.2749*\x});
			\draw[thick, blue, domain=0:0.1100] plot (\x+8.4503,{-2*\x*\x+0.2199*\x});
			\draw[thick, blue, domain=0:0.0880] plot (\x+8.5603,{-2*\x*\x+0.1759*\x});
			\draw[thick, blue, domain=0:0.0704] plot (\x+8.6483,{-2*\x*\x+0.1407*\x});
			\draw[thick,blue] (8.7187,0) -- (9,0);
			\draw[thick, red] (9,0) -- (14,0);
			\draw[thick, green, domain=0:6] plot (\x+14,{1/6*\x*\x-1/54*\x*\x*\x});
			\draw[dashed, thick] (9,-0.5) -- (9,2.5);
			\draw[dashed, thick] (14,-0.5) -- (14,2.5);
	\end{tikzpicture}
	\caption{An illustration of a Zeno control execution.}\label{fig:zeno_control}
\end{figure}
The objective of this work is to study solutions to the following controlled version of \eqref{eq:HDS}
\begin{equation}\label{eq:controlled_HDS}
	\begin{cases}
		\dot{x} = X(x,u), & x\in M\setminus \mathcal{S}, \\
		x^+ = \Delta(x^-), & x\in\mathcal{S},
	\end{cases}
\end{equation}
i.e. the flow is controlled but both the guard and reset are fixed which is a common model for legged locomotion \cite{898695,1166523,bipedal}.

For the hybrid control system, \eqref{eq:controlled_HDS}, we wish to study solutions to the optimal control problem with cost
\begin{equation}\label{eq:cost}
    J(x_0,u(\cdot)) = \int_0^{T_f} \, \ell(x(s),u(s))\, ds,
\end{equation}
subjected to the fixed end-point conditions $x(0) = x_0$ and $x(T_f) = x_f$. A necessary condition for minimizing \eqref{eq:cost} subject to \eqref{eq:controlled_HDS} is the \textit{hybrid maximum principle} \cite{DMITRUK2008964,Pakniyat2017OnTR,4303244}. Unfortunately, this procedure breaks down if and when a Zeno trajectory occurs. However, the hybrid maximum principle dictates that the co-state equation satisfied a ``Hamiltonian jump condition'' and the resulting motion is analogous of an elastic impact system and, as such, Zeno will almost never happen \cite{2111.11645}.

The contribution of this work is to demonstrate that ``almost never Zeno'' does not mean ``no Zeno.'' In the specific case of the controlled bouncing ball with dissipation at impacts, if the terminal time is sufficiently large, it is advantageous to only apply controls around the terminal time. An illustration of this control procedure is shown in Fig. \ref{fig:zeno_control}. 

Preliminaries in both hybrid systems and their accompanying control systems are outlined in \S\ref{sec:prelim}. The classical hybrid maximum principle is presented in \S\ref{sec:hmp}. A detailed examination of this problem applied to the case study of the bouncing ball is performed in \S\ref{sec:bball} where it is explicitly shown that in certain circumstances the Zeno trajectory shown in Fig. \ref{fig:zeno_control} describes the optimal solution. Finally, conclusions and future directions are discussed in \S\ref{sec:conclusions}.
\section{PRELIMINARIES}\label{sec:prelim}
A hybrid dynamical system is composed of a continuous flow generated by a vector field, and a discrete reset map which is activated whenever the flow intersects the guard. The specific definition used throughout is given below.
 
\begin{definition}\label{def:hybrid_system}
	A hybrid dynamical system is a 4-tuple, $\mathcal{H} = (M,\mathcal{S},X,\Delta)$ where
	\begin{enumerate}
		\item $M$ is a (finite-dimensional) manifold,
		\item $\mathcal{S}\subset M$ is an embedded co-dimension 1 submanifold,
		\item $X:M\to TM$ is a vector-field, and
		\item $\Delta:\mathcal{S}\to M$ is a map.
	\end{enumerate}
\end{definition}
Throughout this work all of the data will be assumed to be sufficiently smooth. The manifold $M$ will be referred to as the \textit{state-space}, $\mathcal{S}$ as the \textit{guard}, and $\Delta$ as the \textit{reset}.

The equations of motion for the hybrid system are characterized by their constituent discrete and continuous parts. Off the guard $\mathcal{S}$, the hybrid flow will evolve according to the vector field. Whereas on $\mathcal{S}$ , the hybrid flow instantaneously jumps according to the reset map $\Delta$.
\begin{equation}\label{eq:eom_hybrid}
    	\begin{cases}
			\dot{x}(t) = X(x(t)), & x(t)\not\in \mathcal{S}, \\
			x(t^+) = \Delta(x(t^-)), & x(t^-)\in \mathcal{S}.
		\end{cases}
\end{equation}
The major qualitative feature that is unique to hybrid systems is that of \textit{Zeno}.
\begin{definition}
    Let $\phi_t^{\mathcal{H}}$ be the flow of \eqref{eq:eom_hybrid}. Then a point $x\in M$ has a Zeno trajectory if there exists an increasing sequence of times $\{t_i\}_{i = 1}^{\infty}$ such that $\phi_{t_i}^{\mathcal{H}}(x) \in \mathcal{S}$ for all $i$ and the limit $\lim_{i \to \infty}t_i = t_{\infty}$ exists and is finite.
\end{definition}

In order to provide a theoretical framework for performing optimal control in the context of hybrid systems, we extend the notion of a hybrid dynamical system (Definition \ref{def:hybrid_system}) to include a control variable, denoted by $u$, to the continuous component of the dynamics. No control over either the reset or guard will be assumed.
\begin{definition}
	A hybrid control system is a 5-tuple, $\mathcal{HC}=(M,\mathcal{U},\mathcal{S},X,\Delta)$ where
	\begin{enumerate}
		\item $M$ is a (finite-dimensional) manifold,
		\item $\mathcal{S}\subset M$ is an embedded co-dimension 1 submanifold,
		\item $\mathcal{U}\subset \mathbb{R}^m$ is a closed subset of admissible controls,
		\item $X:M\times V\to TM$ is where $\mathcal{U}\subset V$ is an open neighborhood, and
		\item $\Delta:\mathcal{S}\to M$ is a map.
	\end{enumerate}
\end{definition}
Again, it will be implicitly assumed that all the data are sufficiently smooth.

The optimal control problem of interest will be the following.
\begin{prob}\label{problem}
    For a given a hybrid control system $(M, \mathcal{U}, \mathcal{S}, X, \Delta)$, minimize the cost functional:
	\begin{gather*}
		J:M\times \mathcal{U}^{[0,T_f]}\times[0,T_f]\to\mathbb{R}, \\
		J(x_0,u(\cdot),s) = \int_s^{T_f} \, \ell(x(t),u(t)) \, dt,
	\end{gather*}
	subject to the boundary conditions $x(s)=x_0$ and $x(T_f) = x_f$, along with the controlled dynamics
	\begin{equation*}
		\begin{cases}
			\dot{x}(t) = X(x(t),u(t)), & x(t)\not\in \mathcal{S}, \\
			x(t^+) = \Delta(x(t^-)), & x(t^-)\in \mathcal{S}.
		\end{cases}
	\end{equation*}
\end{prob}

\section{HYBRID MAXIMUM PRINCIPLE}\label{sec:hmp}
As Problem \ref{problem} obeys the principle of optimality, it is amendable to the usual ideas of optimal control theory; there exists both a hybrid Hamilton-Jacobi-Bellman equation \cite{Haddad} and a hybrid maximum principle \cite{DMITRUK2008964,Pakniyat2017OnTR,4303244}.
Although it will only provide necessary conditions, we will use the hybrid maximum principle to study Problem \ref{problem}.
The continuous part of the optimal trajectory follows the flow of the optimal Hamiltonian, just as in the classical, non-hybrid, case. For hybrid control system  $\mathcal{HC}=(M,\mathcal{U},\mathcal{S},X,\Delta)$ with cost functional $J$ and running cost $\ell$, the optimal Hamiltonian, $\hat{H}:T^*M\to\mathbb{R}$, is given by:
\begin{align*}
    \hat{H}(x, p) &= \min_{u \in U}\,\tilde{H}(x, p, u) \\
    &= \min_{u \in U}\, \ell(x, u) + \langle p, X(x, u) \rangle,
\end{align*}
where $\langle p,X\rangle$ denotes the natural pairing between a co-vector and a vector. 
At impacts, the optimal flow is discontinuous. 
As the state variables jump according to the reset map, $\Delta$, the co-states jump according to the \textit{extended reset map}, $\tilde{\Delta}$, which satisfies:
\begin{equation}\label{eq:extended_reset}
	\left(\mathrm{Id}\times\tilde{\Delta}\right)^*\vartheta_{\hat{H}} = \iota^*\vartheta_{\hat{H}},
\end{equation}
such that the following diagram is commutative
\begin{equation*}
	\begin{tikzcd}[ampersand replacement=\&,column sep=large, row sep=large]
		\mathbb{R}\times \mathcal{S}^* \arrow[r, "\mathrm{Id}\times\tilde{\Delta}"] \arrow[d, "\mathrm{Id}\times\pi_M"] \& \mathbb{R}\times T^*M \arrow[d, "\mathrm{Id}\times\pi_M"] \\
		\mathbb{R}\times \mathcal{S} \arrow[r,"\mathrm{Id}\times\Delta"] \& \mathbb{R}\times M
	\end{tikzcd}
\end{equation*}
where $\vartheta_{\hat{H}} = p_i\cdot dx^i - \hat{H} \cdot dt$ is the action form, $\pi_M:T^*M\to M$ is the cotangent projection, and $\iota:\mathcal{S}^*\hookrightarrow T^*M$ is the inclusion map of the \textit{extended guard},
\begin{equation*}
    \mathcal{S}^* = \left\{ (x,p) \in T^*M|_\mathcal{S} : dh_x\left( \frac{\partial H}{\partial p}\right) > 0 \right\}.
\end{equation*}

The hybrid maximum principle states that a necessary condition for a trajectory, $x(t)$, to be a solution to Problem \ref{problem}, it must have the form $x(t) = \pi_M(\gamma(t))$ where $\gamma(t)$ is an integral curve of 
\begin{equation*}
    \begin{cases}
        \dot{x} = \dfrac{\partial \hat{H}}{\partial p}, \quad \dot{p} = -\dfrac{\partial \hat{H}}{\partial x}, & x\not\in \mathcal{S}, \\[2ex]
        (x^+,p^+) = \tilde{\Delta}(x^-,p^-), & x\in \mathcal{S},
    \end{cases}
\end{equation*}
such the boundary conditions are satisfied: $\pi_M(\gamma(0)) = x_0$ and $\pi_M(\gamma(T_f))=x_f$.

Under the assumption that there exists a unique extended reset map that satisfies the conditions above (one does not generally expect this to be true; see \cite{2111.11645} for a more detailed analysis on issue), along with some regularity assumptions, it can be proved that the set of integral curves obeying the hybrid maximum principle that are Zeno have measure zero.
However, this does not imply that the actual optimal trajectory will not be Zeno in any specific example.
Indeed, in the next section we show that the Zeno phenomena is present even in the simple case of a bouncing ball. 
\section{CASE STUDY: THE BOUNCING BALL}\label{sec:bball}
\begin{figure}
	\centering
	\includegraphics[width=0.9\columnwidth]{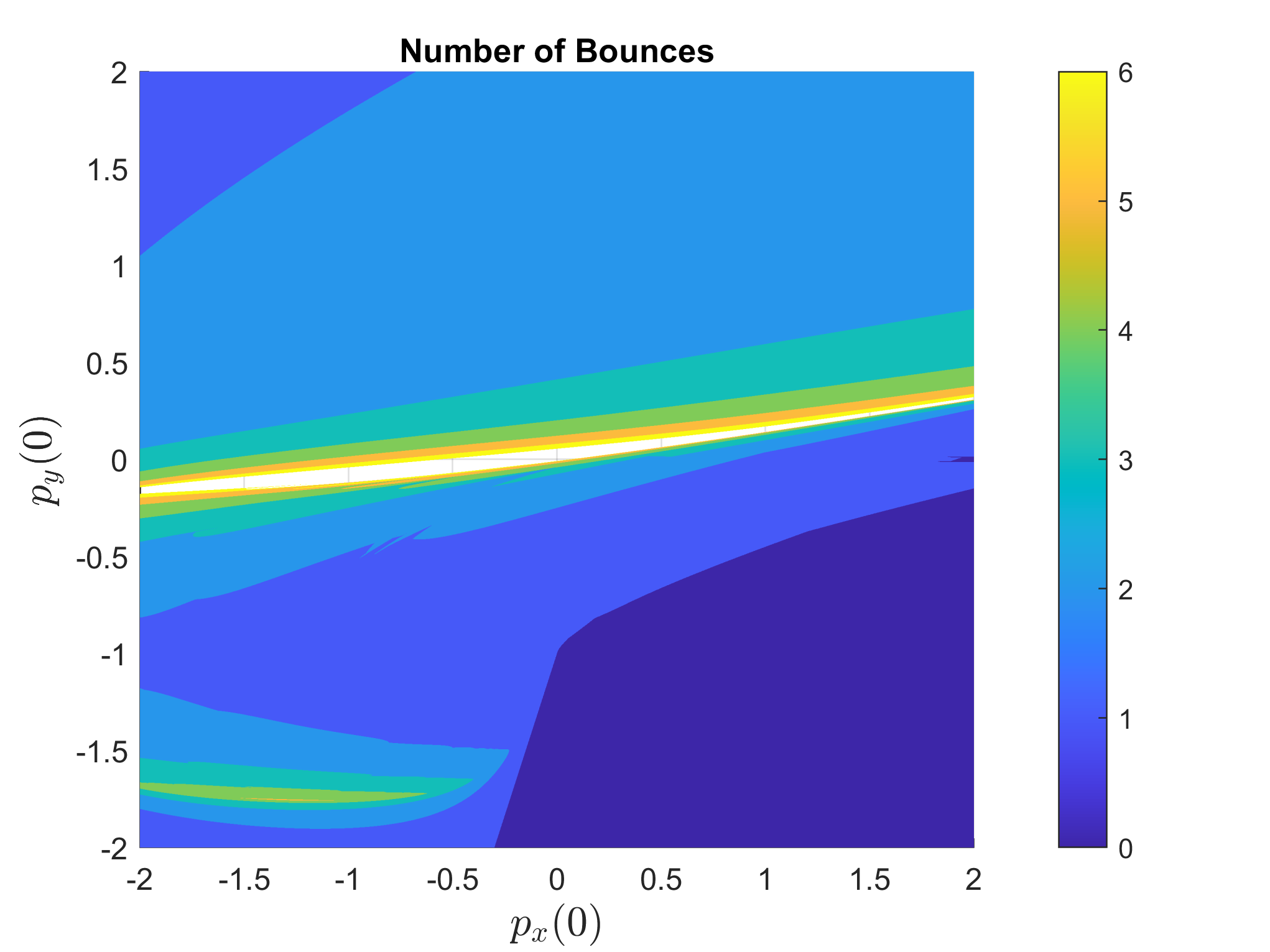}
	\caption{The number of bounces the optimal trajectory undergoes as a function on the initial momentum. The initial condition was $x(0) = 0.1$ and $y(0)=0$. In a neighborhood of the origin, the number of bounces is unbounded and so the image is cropped to only display the regions where there are less than or equal to six bounces. The location where the number of resets grows without bound appears to be a curve with zero volume which is in agreement with the face that almost no trajectories are Zeno.}
	\label{fig:num_bounces}
\end{figure}

\begin{figure}
    \centering
    \includegraphics[width=0.9\columnwidth]{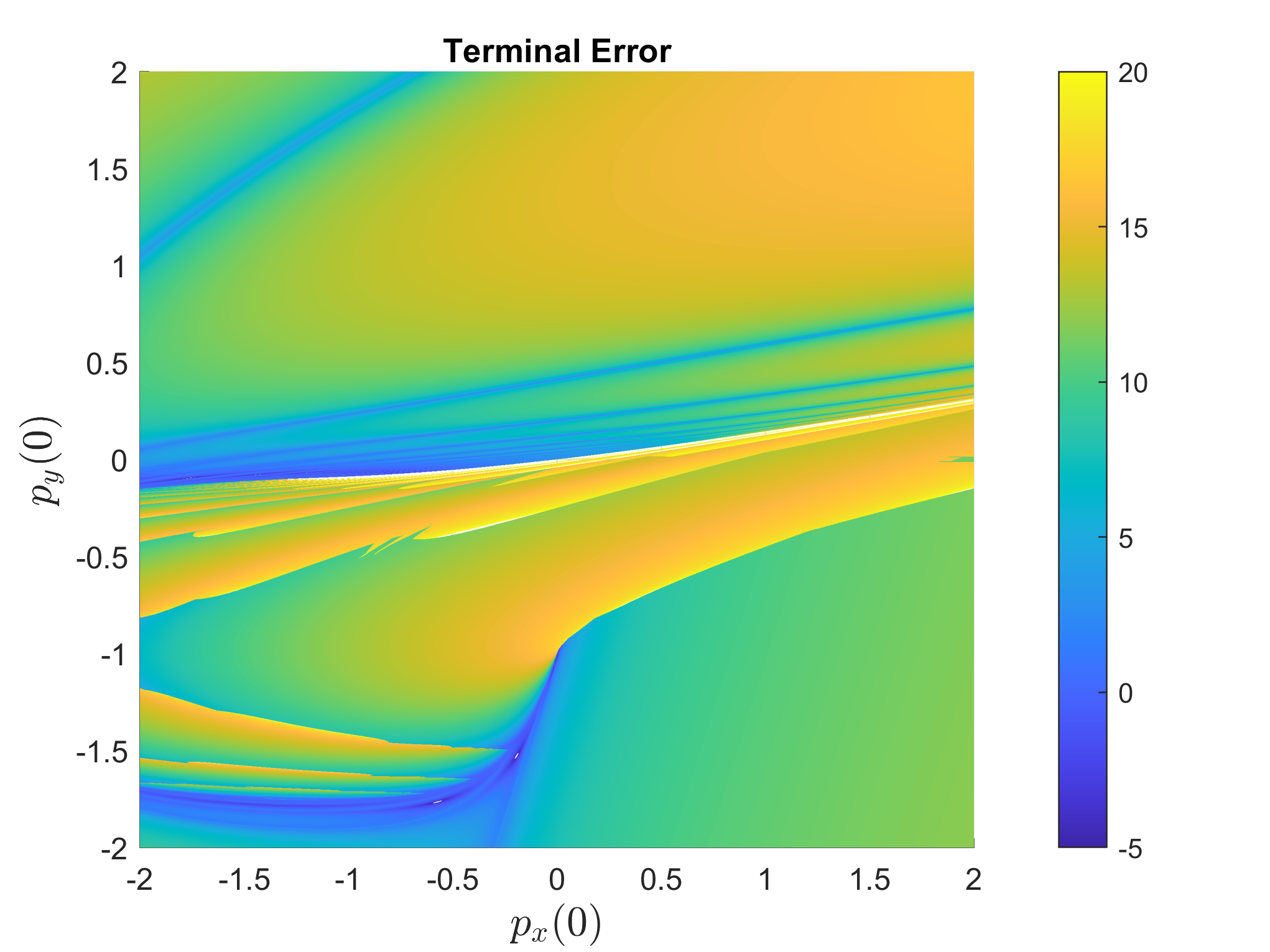}
    \caption{A plot of the terminal error as a function on the initial momentum subject to the same initial conditions as Fig. \ref{fig:num_bounces}. The plot shows the logarithm of the error, $\ln\left( (x(T_f)-x_f)^2+(p(T_f)-p_f)^2\right)$. The error is both discontinuous and spans many orders of magnitude.}
    \label{fig:log_error}
\end{figure}

Consider the case of the bouncing ball with mass $m$ and gravity $g$ subject to impact dissipation. The controlled continuous dynamics are given by
\begin{equation}\label{eq:cont_bball}
	\begin{split}
		\dot{x} &= \frac{1}{m}p, \\
		\dot{p} &= -mg + u,
	\end{split}
\end{equation}
where $x>0$ is the height of the ball, $p$ its momentum, and an admissible controls being any real number, $u\in\mathbb{R}$.
The reset is inelastic with coefficient of restitution $0<c<1$,
\begin{equation}\label{eq:reset_bball}
	\Delta:(x,p) \mapsto (x,-c^2p),
\end{equation}
and occurs on the guard,
\begin{equation*}
	\mathcal{S} = \left\{ (x,p) : x=0, \; p<0\right\}\subset T^*\mathbb{R}.
\end{equation*}
The cost is will be chosen to simply be
\begin{equation}\label{eq:cost_bball}
	J = \int_0^{T_f} \, \frac{1}{2}u^2 \, dt,
\end{equation}
with terminal conditions $x(T_f) = x_f$ and $p(T_f) = p_f$. The specific parameters used throughout this section are shown in Table \ref{tab:parameter_cart}.
\begin{table}
	\centering
	\vspace{0.05in}
	{\renewcommand{\arraystretch}{1.2}
	\begin{tabular}{c|c}
		Parameter& Value  \\ \hline
		$m$ & 1 \\
		$g$ & 1 \\
		$T_f$ & $10$ \\
		$c$ & $0.75$ \\
		$x_f$ & 1 \\
		$p_f$ & 0
	\end{tabular}}
	\caption{Parameters used for the bouncing ball case study.}
	\label{tab:parameter_cart}
\end{table}
The optimal Hamiltonian is
\begin{equation}\label{eq:opt_Ham}
    \begin{split}
        \hat{H} &= \min_{u} \left[ \frac{1}{2}u^2 + \frac{1}{m}p_xp + (-mg+u)p_p\right] \\
        &= -\frac{1}{2}p_p^2 + \frac{1}{m}p_xp - mgp_p,
    \end{split}
\end{equation}
where $p_x$ and $p_p$ are the co-states corresponding to $x$ and $p$ respectively, and the optimal control is $u^* = -p_p$.
Therefore, in-between impacts, the continuous dynamics evolve according to
\begin{gather}\label{eq:bball_flow}
    \dot{x} = \frac{1}{m}p, \quad \dot{p}_x = 0, \\
    \dot{p} = -mg-p_p, \quad \dot{p}_p = -\frac{1}{m}p_x. \nonumber
\end{gather}
The extended reset map is uniquely determined by \eqref{eq:extended_reset}
\begin{equation}\label{eq:bball_reset}
	\begin{split}
		x&\mapsto x, \\
		p&\mapsto -c^2p, \\
		p_x &\mapsto -\frac{1}{c^2}p_x + \frac{m}{2c^2}\frac{p_p^2}{p}\left(1-c^{-4}\right) \\
		&\qquad\qquad\qquad + \frac{m^2g}{c^2}\frac{p_p}{p}\left(1+c^{-2}\right), \\
		p_p &\mapsto -\frac{1}{c^2}p_p.
	\end{split}
\end{equation}
As long as Zeno is not present, the optimal trajectory can be found by concatenating the integrated flow of \eqref{eq:bball_flow} with the gluing condition \eqref{eq:bball_reset}. Numerically, this was implemented via a na\"{i}ve shooting method to ensure that the terminal conditions are satisfied. The resulting shooting problem is highly discontinuous, see Figures \ref{fig:num_bounces} and \ref{fig:log_error}.

\begin{figure}
	\centering
	\includegraphics[width=0.9\columnwidth]{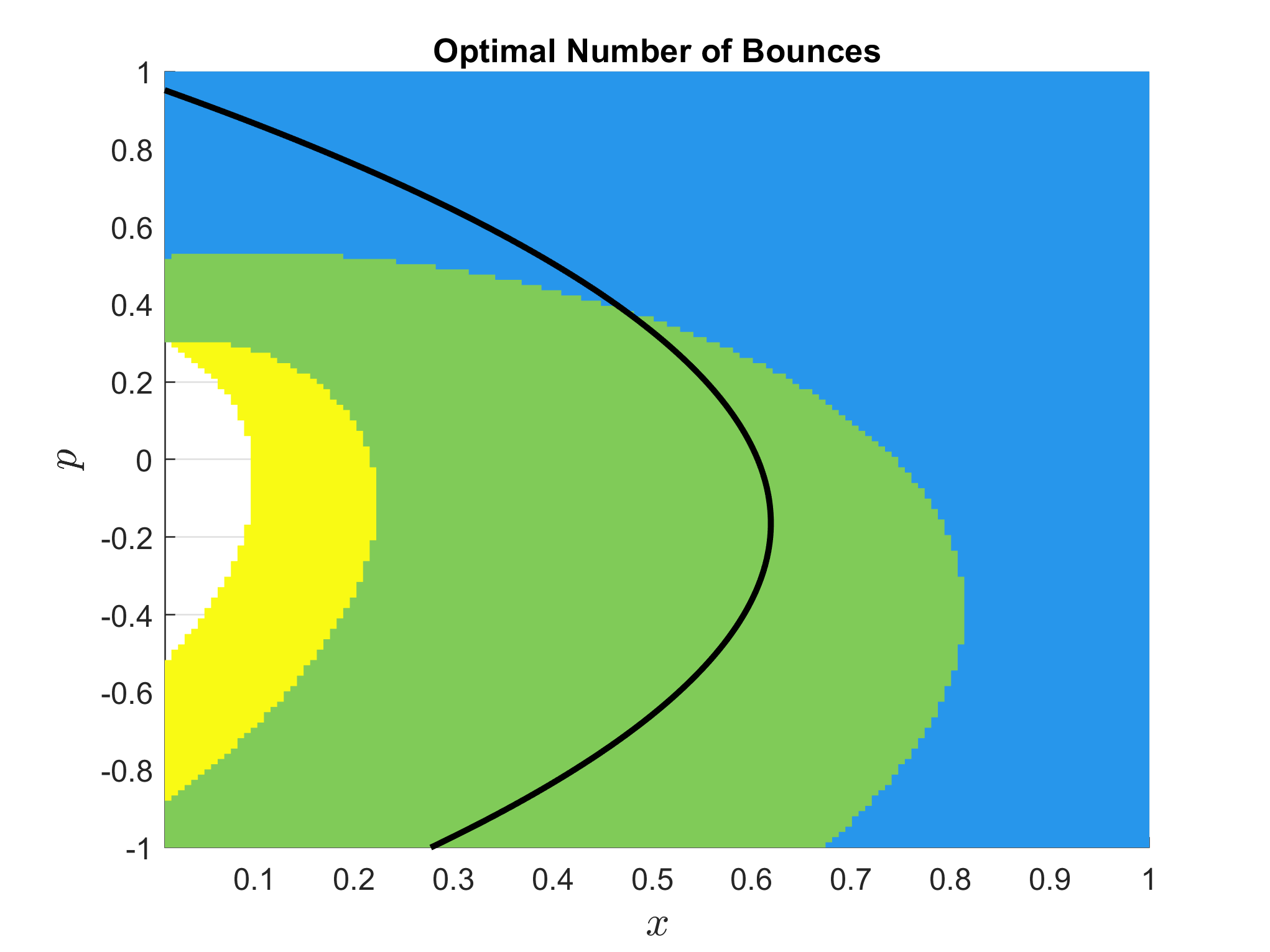}
	\caption{The black curve denotes the boundary of where the Zeno control scheme becomes locally optimal, i.e. where $T_f = \zeta(x,p) + \sqrt{6/g}$. The colored regions provide the optimal number of bounces; blue is three, green is four, and yellow is five. The white region is where the Zeno control is optimal.}
	\label{fig:num_bounces_opt}
\end{figure}

The trajectories are numerically solved by integrating \eqref{eq:bball_flow} via Matlab's \texttt{ode45} adaptive Runge-Kutta solver and event detection for applying the reset map \eqref{eq:bball_reset}. The shooting problem is solved by calling Matlab's \texttt{fsolve}, a trust-region dogleg algorithm, with $625=25^2$ initial seeds chosen uniformly, $-2\leq p_x(0),p_y(0)\leq 2$. The colored part of Fig. \ref{fig:num_bounces_opt} and the corresponding piece of Fig. \ref{fig:value_fcn} report the found optimal number of bounces and cost, respectively.

\subsection{The Zeno Control Scheme}
The above synthesis breaks down when the trajectory is Zeno as this occurs when $(x,p)=(0,0)$ and \eqref{eq:bball_reset} fails to be defined. As such, in order to describe the Zeno control scheme as shown inf Fig. \ref{fig:zeno_control}, we must proceed independently of the above procedure.
With the controls turned off, any initial condition of the system (\ref{eq:cont_bball},\ref{eq:reset_bball}) with $u=0$ results in a Zeno trajectory. The time at which the trajectories collapse to Zeno can be calculated in closed-form as the underlying dynamics can be integrated exactly:
\begin{equation*}
    \zeta(x_0,p_0) = \frac{1}{mg}p_0 + \frac{3}{g(1-c^2)}\sqrt{\frac{p_0^2}{m^2} + 2gx_0}.
\end{equation*}
If the controls are turned on $T$ time before termination, i.e. at $t = T_f-T$, the cost to transition between the Zeno state, $(x(t),p(t))=(0,0)$, and the terminal condition, $(x(T_f),p(T_f))=(1,0)$, is exactly
\begin{equation*}
    J(T) = \frac{1}{2}g^2T + 6T^{-3},
\end{equation*}
as the continuous trajectories, \eqref{eq:bball_flow}, can be solved exactly. This cost is minimized when $T^* = \sqrt{6/g}$ to
\begin{equation}\label{eq:Zeno_cost}
    J(T^*) = \frac{2}{3}g\sqrt{6g}.
\end{equation}
Combining all of the above, if $T_f \geq \zeta(x_0,p_0) + \sqrt{6/g}$, then the Zeno control scheme is locally optimal with cost given by \eqref{eq:Zeno_cost}. 

\begin{figure}
	\centering
	\includegraphics[width=0.9\columnwidth]{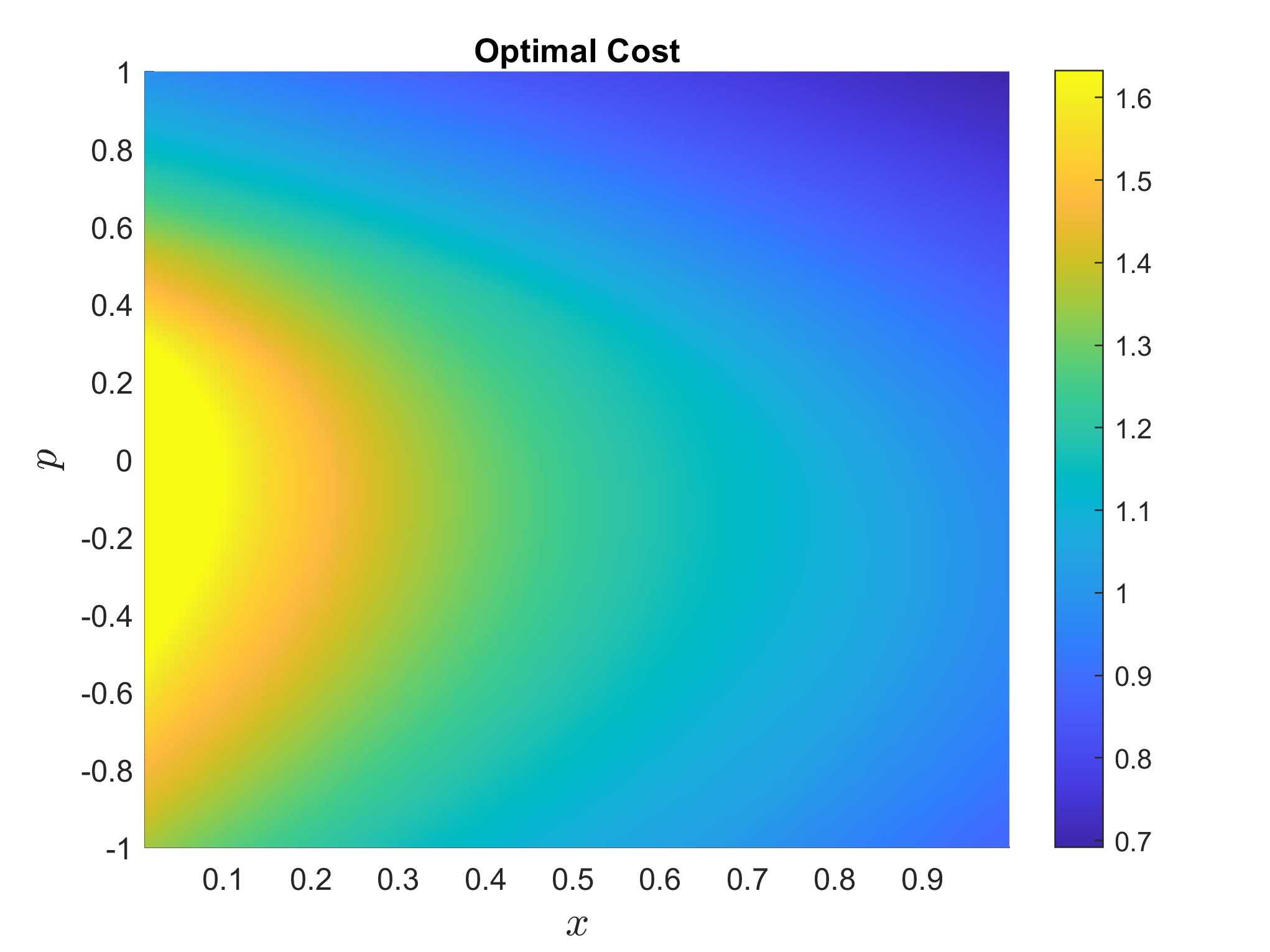}
	\caption{The value function for the bouncing ball by taking the minimum of both the Zeno control cost (where defined) and the cost numerically found via shooting.}
	\label{fig:value_fcn}
\end{figure}

Let $J_{shoot}$ be the optimal cost found by solving the Zeno-free trajectory (\ref{eq:bball_flow},\ref{eq:bball_reset}). Then the true value function is given by
\begin{equation*}
    J_{true} = \min\left\{ J_{shoot}, \frac{2}{3}g\sqrt{6g}\right\}.
\end{equation*}
If $J_{true} \ne J_{shoot}$, then the Zeno trajectory is the optimal solution and the hybrid maximum principle is not a sufficient amount of analysis. The optimal number of bounces is numerically computed in Fig. \ref{fig:num_bounces} while the true value function is in Fig. \ref{fig:value_fcn} where it can be seen that the Zeno trajectory is optimal in a neighborhood around the origin.
\subsection{Numerical Observations}
There are a number of interesting features of Figures \ref{fig:num_bounces_opt} and \ref{fig:value_fcn}. Arguably the most important observation is that the region where the Zeno control scheme is optimal is non-empty. This demonstrates that the hybrid maximum principle is incomplete for solving the optimal control problem.

The other qualitatively interesting observation is the shape of the regions in Fig. \ref{fig:num_bounces_opt}. Despite Figures \ref{fig:num_bounces} and \ref{fig:log_error} being highly complicated and intricate, the regions dictating how many bounces are optimal seem to be separated by smooth curves. Moreover, the region where Zeno is actually optimal is much smaller than the region where it is only locally optimal.

There are quantitative features of Fig. \ref{fig:value_fcn} that are also important to notice. First of all, the value function has the constant value from \eqref{eq:Zeno_cost} when ever Zeno is optimal; this means that the Zeno region is Fig. \ref{fig:value_fcn} is constant in both space and time. Constant solutions are a valid solution for the corresponding Hamilton-Jacobi-Bellman equation,
\begin{equation*}
    \frac{\partial J}{\partial t} - \frac{1}{2}\left(\frac{\partial J}{\partial p}\right)^2 + \frac{1}{m}p\frac{\partial J}{\partial x} - mg\frac{\partial J}{\partial p} = 0,
\end{equation*}
and the accompanying boundary conditions corresponding to the extended reset \eqref{eq:bball_reset}.
On top of this, even though the extended reset equations do not properly handle the Zeno trajectory, the optimal Hamiltonian \eqref{eq:opt_Ham} remains preserved. As controls are initially turned off, we have $p_x(0)=p_y(0)=0$ and thus $\hat{H} = 0$. When controls are turned on at $t^* = T_f - \sqrt{6/g}$, the co-states are
\begin{equation*}
    p_x(t^*) = -\sqrt{2/3}g^{3/2}, \quad p_y(t^*) = -2g,
\end{equation*}
and the Hamiltonian is
\begin{equation*}
    \begin{split}
        \hat{H} &= -\frac{1}{2}p_y^2 - gp_y = -\frac{1}{2}4g^2 + 2g^2 = 0.
    \end{split}
\end{equation*}
\section{CONCLUSIONS}\label{sec:conclusions}
\begin{figure}
	\centering
	\includegraphics[width=0.9\columnwidth]{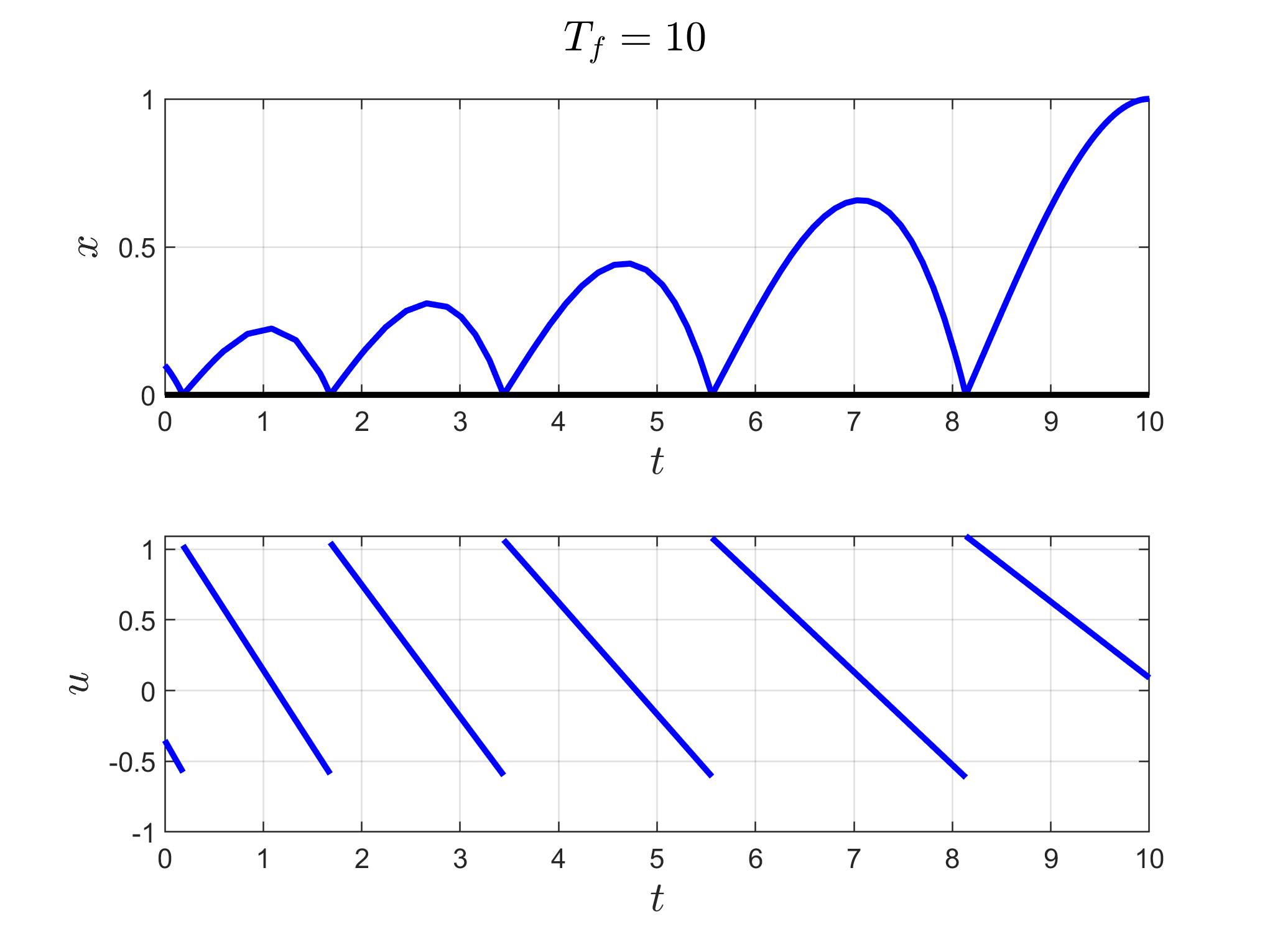}
	\caption{The optimal controls for an initial condition with five bounces as dictated by the hybrid maximum principle.}
	\label{fig:five_path}
\end{figure}

\begin{figure}
	\centering
	\includegraphics[width=0.9\columnwidth]{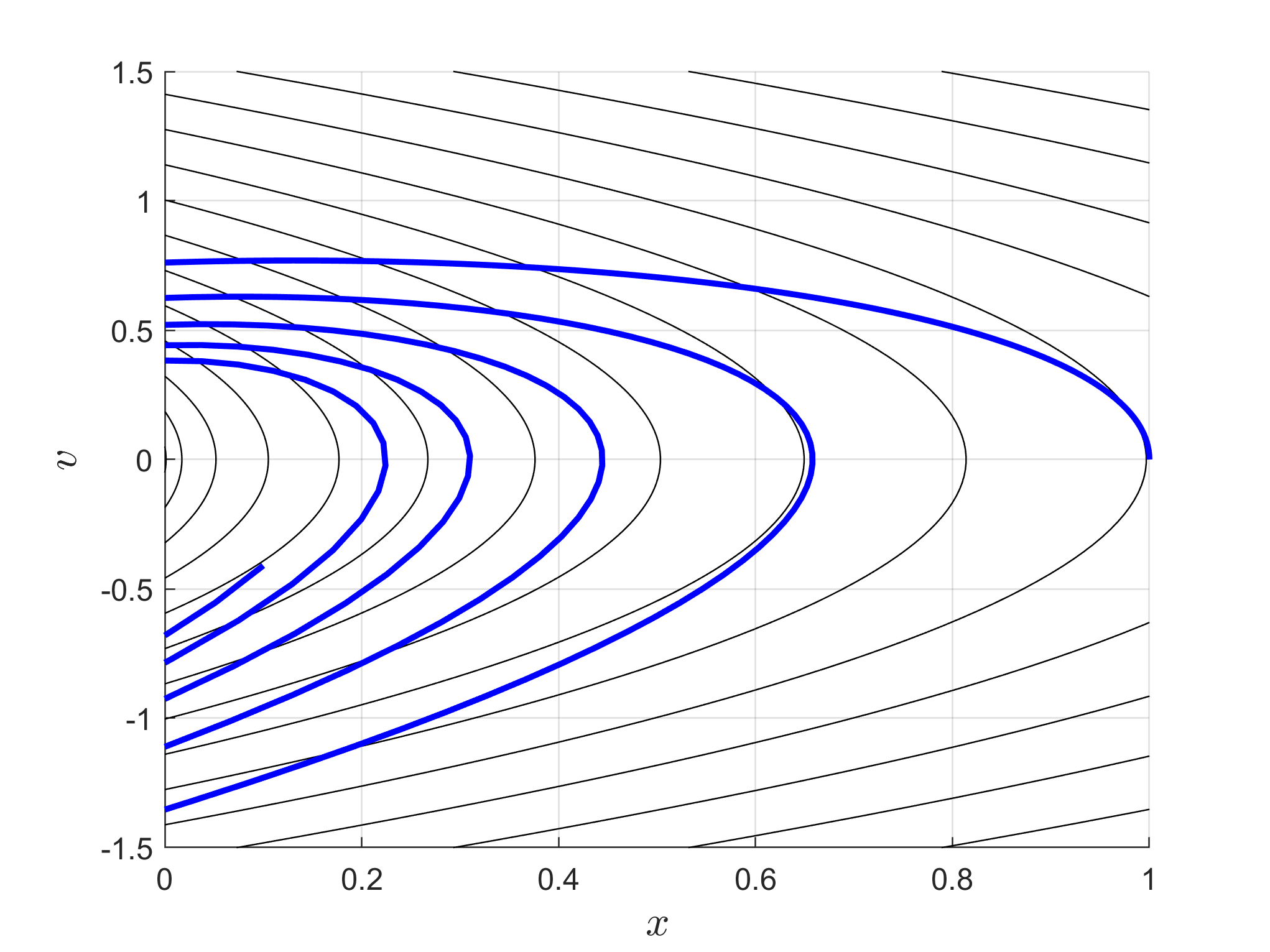}
	\caption{The phase portrait corresponding to Fig. \ref{fig:five_path}.}
	\label{fig:five_phase}
\end{figure}

\begin{figure}
	\centering
	\includegraphics[width=0.9\columnwidth]{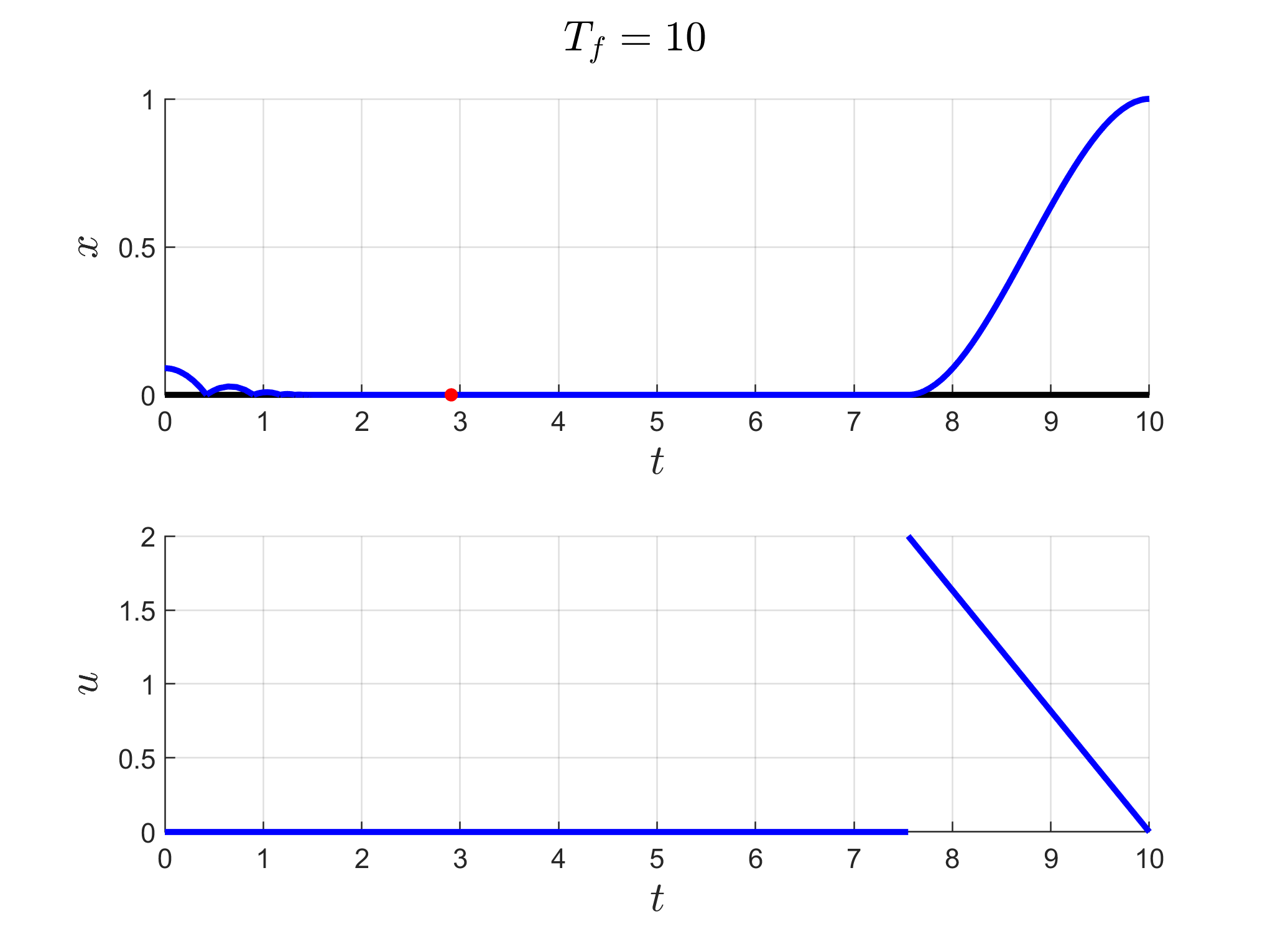}
	\caption{The optimal controls for an initial conditions with a Zeno trajectory. This image displays the same features as Fig. \ref{fig:zeno_control}.}
	\label{fig:zeno_path}
\end{figure}

\begin{figure}
	\centering
	\includegraphics[width=0.9\columnwidth]{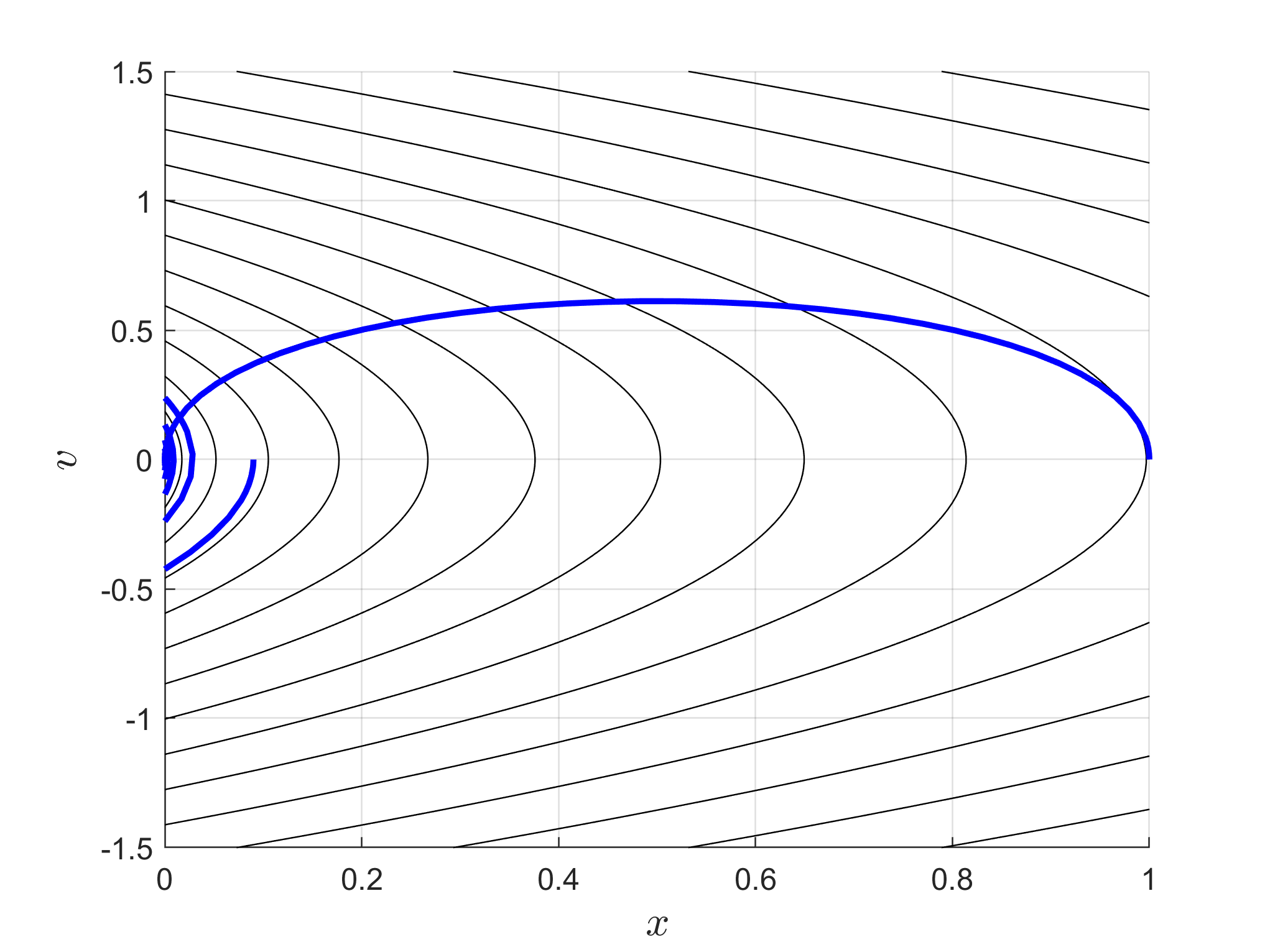}
	\caption{The phase portrait corresponding to Fig. \ref{fig:zeno_path}.}
	\label{fig:zeno_phase}
\end{figure}

This work studied the existence of Zeno trajectories in optimally controlled hybrid systems by studying the pedagogical example of the dampened bouncing ball. When the terminal time was sufficiently large and the initial conditions sufficiently close to the Zeno state, it turned out to be advantageous to allow the trajectory to be Zeno and only apply the controls towards the end; a sample trajectory is shown in Figures \ref{fig:zeno_path} and \ref{fig:zeno_phase} which is qualitatively distinct from the solutions arising from the hybrid maximum principle as shown in Figures \ref{fig:five_path} and \ref{fig:five_phase}.


Although the bouncing ball example provides a valuable insight into the mechanisms of Zeno trajectories and their importance in solving the optimal control problem, a more general framework for dealing with Zeno trajectories still needs to be developed and applied to more complicated examples. In particular, note that for the bouncing ball case, the hybrid flow can be analytically computed, and there is a closed form solution for the time at which trajectories become Zeno. However, in general, we might not be able to explicitly compute either of these. As such it would be desirable to design an algorithm that would output whether a Zeno trajectory can be a solution to the optimal control problem, whether it is unique, and compute the time to Zeno if this is the case. 

Finally, our numerical experiments (see Fig. \ref{fig:num_bounces_opt}) indicate that the interface between the areas with different number of bounces are smooth curves. It would be of great interest to find an analytical formula, or merely approximations, to characterize these curves. Additionally, the region in which the Zeno control scheme becomes locally optimal is larger than the area where Zeno is the globally optimal solution to the control problem. This is not a surprise as the maximum principle only gives necessary conditions for optimality and not sufficient. Still, it would be desirable to find out what causes this mismatch, and provide sufficient conditions for Zeno.


\addtolength{\textheight}{-12cm}   



\section*{ACKNOWLEDGMENT}
We thank Prof. Alexander Vladimirsky for his interest and encouragement of this work.

\bibliographystyle{ieeetr}
\bibliography{references.bib}

\end{document}